\documentclass[12pt]{article}
\usepackage{amsmath,amssymb}
\usepackage{geometry}
\usepackage{amsmath, amsfonts,amsthm,amssymb,enumerate,amscd}

\usepackage{times}

\usepackage[active]{srcltx}
\usepackage[utf8]{inputenc}
\usepackage{algorithm}
\usepackage{stmaryrd}
\usepackage{array}
\usepackage{algorithmicx}
\usepackage{algpseudocode}
\usepackage{graphicx}
\usepackage{letltxmacro}

\usepackage{color}
\usepackage[table,xcdraw,svgnames,dvipsnames]{xcolor}
\definecolor{light-salmon}{RGB}{255,140,120}
\usepackage[colorlinks]{hyperref}
\hypersetup{
	allbordercolors=.,
	citecolor=DarkGreen,
	linkcolor=DarkRed,
	urlcolor=NavyBlue
}

\graphicspath{{./pics/}{./metapost/}}

\usepackage{enumitem}

\numberwithin{equation}{section}
\theoremstyle{plain}
\newtheorem{thm}{Theorem}

\newtheorem{prop}[thm]{Proposition}
\newtheorem{deff}[thm]{Definition}
\theoremstyle{definition}
\newtheorem{rem}[thm]{Remark}
\newtheorem{conj}[thm]{Conjecture}

\newcommand{\Per}{\operatorname{Per}}

\renewcommand{\Bbb}{\mathbb}
\newcommand{\bo}[1]{{\bf #1}}

\newcommand{\F}{\mathcal{F}}
\newcommand{\A}{\mathcal{A}}
\newcommand{\C}{\mathcal{C}}
\newcommand{\trc}{{three-$r$-cap}}

\makeatletter
\DeclareFontFamily{U}{tipa}{}
\DeclareFontShape{U}{tipa}{m}{n}{<->tipa10}{}
\newcommand{\arc@char}{{\usefont{U}{tipa}{m}{n}\symbol{62}}}%

\newcommand{\arc}[1]{\mathpalette\arc@arc{#1}}

\newcommand{\arc@arc}[2]{%
	\sbox0{$\m@th#1#2$}%
	\vbox{
		\hbox{\resizebox{\wd0}{\height}{\arc@char}}
		\nointerlineskip
		\box0
	}%
}
\makeatother

\title{Shape optimization under width constraint: the Cheeger constant and the torsional rigidity}

\author{Beniamin Bogosel}

\begin{document}
	
\maketitle

\begin{abstract}
In this article it is shown that the equilateral triangle maximizes the Cheeger constant and minimizes the torsional rigidity among shapes having a fixed minimal width. The proof techniques use direct comparisons with simpler shapes, consisting of disks with three disjoint caps. Comparison results for harmonic functions help establish that in non-equilateral configurations the shape derivative has an appropriate sign, contradicting optimality. 
\end{abstract}



\section{Introduction}

Shape optimization problems deal with minimization of functionals depending on the geometry of a shape subjected to various constraints. The canonical formulation is
\begin{equation}\label{eq:sh-opt}
	 \min_{\Omega \in \mathcal A}\F(\Omega).
\end{equation}
The set $\mathcal A$ of admissible shapes includes the constraints imposed on the geometry of $\Omega$. Depending on the continuity properties of $\F$ and on the compactness of the family $\A$, the existence of solutions for problem \eqref{eq:sh-opt} is more or less straightforward. 

Recently, many works in the literature dealt with problems related to convexity, diameter or width constraints. In this case, compactness properties related to the Hausdorff convergence of convex shapes imply existence of solutions for a wide range of problems. Details and references can be found in \cite[Section 2]{henrot-pierre-english} and \cite[Section 3]{AntunesBogosel22}. While the existence of solutions can be, in general, obtained using the classical method of calculus of variations, qualitative properties and identification of solutions is often difficult, since optimality conditions are not straightforward to exploit \cite{LNconv09}, \cite{LNPconv12}. To give intuition regarding optimization problems under convexity and width constraints numerical methods were proposed in \cite{AntunesBogosel22} and \cite{Bogosel-Convex}. The simulations in \cite{Bogosel-Convex} show for example that it is likely that the eigenvalues $\lambda_k(\Omega)$, solving $-\Delta u = \lambda_k(\Omega) u$ with $u \in H_0^1(\Omega)$ are maximized for:
\begin{itemize}[noitemsep]
	\item the Reuleaux triangle in the class of shapes with fixed constant width
	\item the equilateral triangle in the class of shapes with fixed minimal width.
\end{itemize}

The objective of this article is to bring evidence towards the validity of the second claim mentioned above. To fix notations, denote by $K$ a convex domain in $\Bbb{R}^2$, which by convention will be closed and bounded.  Let us define the width of $K$ by $w : [0,2\pi] \to \Bbb{R}$, $w(\theta)$ being equal to the distance between two supporting lines to $K$ orthogonal to $(\cos \theta, \sin \theta)$. A shape $K$ has minimal width $1$ if for every $\theta \in [0,2\pi]$ we have $w(\theta) \geq 1$. Furthermore, a shape has width $w_0>0$ if
\[ w_0 = \min_{\theta \in [0,2\pi]} w(\theta),\]
equivalently, the minimum of the width function for $K$ is equal to $w_0$. 

A classical shape optimization result under minimal width constraint is attributed to Pal \cite{Pal1921} and states that:

\begin{thm}\label{thm:pal}
	The equilateral triangle of unit minimal width minimizes the area among all planar convex shapes having unit minimal width. 
\end{thm}

A proof of this result by direct comparison with a simpler shape is provided in \cite[Exercise 6-4]{yaglom-boltjanskii}. Indeed, in the same reference, Exercise 6-2 it is shown that the inradius of shapes with minimal width $1$ verify $r \in [1/3,1/2]$. For a shape $K$ of unit width having inradius $r<1/2$ there exists a particular simpler shape $T_{ABC}$ contained in $K$: the convex hull of a disk of radius $r$ with center $O$ and three points $A,B,C$ situated at distance $1-r$ from $O$. Since $T_{ABC} \subset K$ we have the obvious inequality $|T_{ABC}|\leq |K|$. Then a simple one dimensional analysis with respect to $r$ shows that the equilateral triangle has minimal area. Recently, the quantitative version of the Pal inequality was investigated in \cite{Pal-three-versions} using multiple points of view, including comparisons with the set $T_{ABC}$ recalled above. 

This method solves the problem stated in Theorem \ref{thm:pal} by \bo{direct comparison}. In the following, we show how the same method can be applied to solve the analogue problem for the Cheeger constant and the torsional rigidity.  This brings new evidence towards the validity of a similar statement for the Dirichlet Laplace eigenvalues.  

\begin{deff}\label{def:Cheeger}
	Let $K\subset \Bbb{R}^2$ be a convex domain with non void interior. The associated Cheeger constant is given by 
	\[ h(K) = \inf \left\{ \frac{\Per(E)}{|E|}: E \text{ convex }, E \subset K, |E|>0.\right\}\]
	Any set realizing the infimum above is called a \emph{Cheeger set} of $K$.
\end{deff}
Some authors consider a more general setting where the infimum is taken over all measurable sets. Nevertheless, in dimension two taking convex hull decreases the ratio $\Per(E)/|E|$, therefore we can restrict the analysis to convex subsets of $K$. 

The definition implies that the Cheeger constant is decreasing with respect to set inclusion: $K \subset K' \Longrightarrow h(K) \geq h(K')$. This suggests that direct comparison methods might work for studying the maximization of the Cheeger constant under minimal width constraint.  In the recent paper \cite{Cheeger-minw} the authors show that indeed the equilateral triangle maximizes the Cheeger constant under width constraint. The result is recalled below and a quick proof using direct comparison will be shown in the following section.

\begin{thm}\label{thm:Cheeger-minw}
	The equilateral triangle of unit width maximizes the Cheeger constant among shapes having unit minimal width.
\end{thm} 

In the same context of optimization under width constraints, in the paper \cite{henrot-lucardesi} the authors show that the Reuleaux triangle maximizes the Cheeger constant in the class of shapes having constant width. An alternate direct proof of this result is given in \cite{bogosel-cheeger}.

Dirichlet Laplace eigenvalues $\lambda_k(\Omega)$ and eigenfunctions $u_k$ solve the problem 
\[ -\Delta u_k = \lambda_k(\Omega) u_k, u_k \in H_0^1(\Omega).\]
More precisely, Dirichlet boundary condition $u_k =0$ holds on $\partial \Omega$.  
For convex domains it is known that the spectrum of the Dirichlet-Laplacian consists only of eigenvalues 
\[ 0<\lambda_1 (\Omega) < \lambda_2(\Omega) \leq ... \to \infty.\]
The simplicity of the first eigenvalue follows from the simple connectedness of convex domains. We choose the first eigenfunction $u_1$ to be positive in the interior of $\Omega$. We refer to \cite[Chapter 1]{henroteigs} for more details. In particular, the characterization of eigenvalues using Rayleigh quotients shows that they are decreasing with respect to set inclusion
\[ K \subset K' \Leftrightarrow \lambda_k(K) \geq \lambda_k(K'),\]
for all $k \geq 1$. Motivated by the numerical simulations in \cite{Bogosel-Convex} and the result of Theorem \ref{thm:Cheeger-minw} the following conjecture is natural. 

\begin{conj}\label{conj:lamk-minw}
	The equilateral triangle of unit width maximizes the Dirichlet-Laplace eigenvalues among shapes having unit minimal width.
\end{conj}

Another shape functional based on partial differential equations is the Dirichlet energy or the torsional rigidity. Let $\Omega$ be a convex shape and consider the PDE
\begin{equation}\label{eq:torsion}
	 -\Delta u = 1, \ w \in H_0^1(\Omega).
\end{equation}
The torsional rigidity is given by
\[ T(\Omega) = \int_\Omega |\nabla u|^2 dx=\int_\Omega u dx.\]
The torsional rigidity behaves in an opposite way comparing to the Dirichlet-Laplace eigenvalues. For example, it is not difficult to notice the following:
\begin{itemize}[noitemsep]
	\item The solution $u$ of \eqref{eq:torsion} is positive in the interior of $\Omega$;
	\item If $\Omega \subset \Omega'$ then $u_\Omega \leq u_{\Omega'}$.
	\item If $\Omega \subset \Omega'$ then $T(\Omega) \leq T(\Omega')$, that is, the torsion is increasing with respect to set inclusion.
\end{itemize}
For more properties and discussion concerning the relation between $T(\Omega)$ and the first eigenvalue $\lambda_1(\Omega)$ see \cite{Henrot2018},  \cite{vandenBerg2020}. When $\Omega$ is convex, the square root $u^{1/2}$ of the torsion fonction solving \eqref{eq:torsion} is concave \cite{Kennington1985}, even strictly concave \cite{Korevaar1987}. 

Following the comparison with classical optimization problems for Dirichlet-Laplace eigenvalues and torsional rigidity and numerical simulations using tools from \cite{Bogosel-Convex}, we formulate the following theorem which is the main result of this paper.

\begin{thm}\label{conj:torsion-minw}
	The equilateral triangle of unit width minimizes the torsional rigidity among shapes having unit minimal width.
\end{thm}

The paper is structured as follows. In Section \ref{sec:Cheeger} the definition of the simple subset $T_{ABC}$ contained in every shape $K$ of minimal unit width will be recalled. This will lead to a straightforward proof for Theorem \ref{thm:Cheeger-minw}. Section \ref{sec:lamk} will present qualitative results regarding the torsional rigidity of three cap sets $T_{ABC}$ which will lead to a proof of Conjecture \ref{conj:torsion-minw}. 

\section{Maximizing the Cheeger constant among shapes with minimal width}
\label{sec:Cheeger}


The ideas presented in in the following are geometric in nature and are inspired from Chapter 6 in the reference \cite{yaglom-boltjanskii}. For simplicity assume that the considered convex shapes $K$ have unit minimal width in the following. It was underlined in the introduction that the inradius of a shape $K$ with unit width verifies $r \in [1/3,1/2]$. If $r = 1/2$ then obviously $K$ will contain the inscribed disk having radius $1/2$. 

Assume that $r<1/2$, implying that the incircle touches $\partial K$ in at least three points $X,Y,Z$. Otherwise it must touch $\partial K$ in two antipodal points saturating the width, leading to $r = 1/2$. Denote by $D_r$ the incircle and by $O$ its center. The difference $K \setminus D_r$ can be partitioned into three regions $R_X, R_Y, R_Z$ such that $X \notin R_X, Y \notin R_Y, Z \notin R_Z$. 

The supporting line $\ell_X$ at $X$ tangent to $D_r$ has an opposite parallel supporting line $\ell_X'$. Since $K$ has minimal width $1$, the distance from $\ell_X'$ to $O$ is at least $1-r$. Thus, there exists a point $A' \in R_X$ such that $OA'\geq 1-r$. Convexity of $K$ implies that there also exists a point $A \in R_X$ such that $OA = 1-r$. Points $B \in R_Y, C\in R_Z$ are defined similarly at distance $1-r$ from $O$. 

\begin{figure}
	\centering 
	\includegraphics[height=0.45\textwidth]{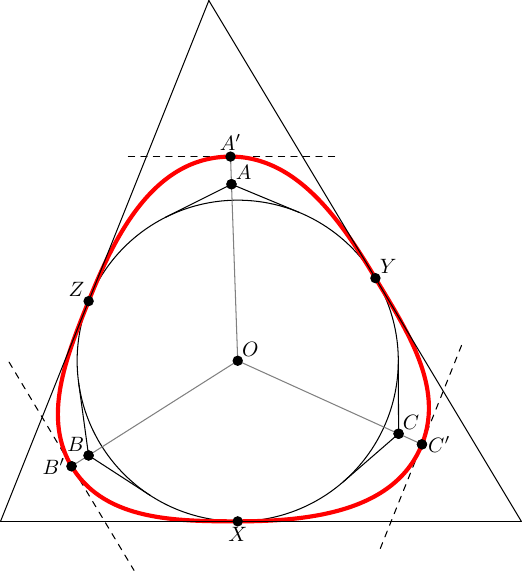}\quad 
	\includegraphics[height=0.45\textwidth]{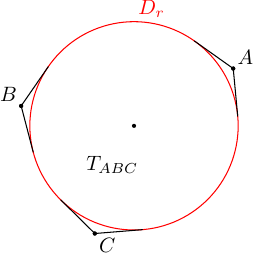}
	\caption{(left) Existence of a three-$r$-cap subset for a shape with given minimal width. (right) Example of a three-$r$-cap shape: convex hull of a disk $D_r$ of center $O$ and radius $r \in [1/3,1/2]$ and points $A,B,C$ such that $AO=BO=CO=1-r$.}
	\label{fig:three-cap}
\end{figure}

Thus, for any shape $K$ of minimal width $1$ and inradius $r<1/2$ there exists a subset $T_{ABC}$ consisting of the convex hull of a disk $D_r$ of radius $R$ with center $O$ and three points $A,B,C$ situated at distance $1-r$ of $O$. Each vertex $A,B,C$ determines a cap in $T_{ABC}$, bounded by the respective vertex and tangents from it to $D_r$. By construction, the three caps from $A,B,C$ are disjoint. For simplicity, call such a set a \bo{three-$r$-cap}. This definition can also be extended to $r=1/2$ noticing that $T_{ABC}$ is simply $D_r$ in this case. See Figure \ref{fig:three-cap} for an illustration.

The monotonicity of the Cheeger constant and Dirichlet eigenvalues with respect to set inclusion imply the following proposition. Denote by $T_e$ the equilateral triangle of unit width (height). 

\begin{prop}\label{prop:reduction-caps}
	(a) Assume that for all \trc{s} $T_{ABC}$ we have $h(T_{ABC}) \leq h(T_e)$. Then the equilateral triangle maximizes the Cheeger constant among shapes of minimal width equal to one.
	
	(b) Assume that for all \trc{s} $T_{ABC}$ we have $\lambda_1(T_{ABC})\leq \lambda_1(T_e)$. Then $\lambda_1$ is maximized by the equilateral triangle among shapes with fixed minimal width.
	
	(c) Assume that for all \trc{s} $T_{ABC}$ we have $T(T_{ABC}) \geq T(T_e)$. Then the torsional rigidity is minimized by the equilateral triangle among shapes with fixed minimal width.
\end{prop}

\emph{Proof:} The proof is immediate from the monotonicity of $h$, $\lambda_k$ and the torsion $T$, coupled with the existence of subsets $T_{ABC}$ for any shape with minimal width equal to one. \hfill $\square$

A classical result of Kawohl and Lachand-Robert states that for convex domains $K$ the Cheeger set can be found by investigating the area of inner parallel sets $K_{-t}$ defined as points in $K$ at distance at least $t$ from $\partial K$. More precisely, $h(K) = 1/t$ where $t>0$ is the solution of the equation $|K_{-t}| = \pi t^2$. It is not difficult to see that inner parallel set of a \trc{} shape is homothetic to itself. Therefore one needs to solve the equation
\[ \frac{(r-t)^2}{r^2} {|T_{ABC}|}=\pi t^2,\]
leading to 
\[ h(T_{ABC}) =  \frac{1}{t} = \frac{1}{r}+\sqrt{\frac{\pi}{|T_{ABC}|}}. \]

Let us compute the area of $T_{ABC}$ in terms of $r$. Denote by $\C_A$ the cap from $A$, the region determined by the tangents from $A$ to $D_r$ and $D_r$. In the same ways define $\C_B, \C_C$ the caps from $B$ and $C$. These caps do not intersect and of course 
\[ |T_{ABC}| = \pi r^2+ |\C_A|+|\C_B|+|\C_C|.\]
Let $A_1, A_2$ be the points where tangents from $A$ intersect $D_r$. Let $\alpha$ be half of the angle $A_1O A_2$, implying $\cos \alpha = \frac{r}{1-r}$, such that $\alpha(r) = \arccos(\frac{r}{1-r})$. Since the inradius $r$ verifies $r \in [1/3,1/2]$ (see \cite{yaglom-boltjanskii} Chapter 6), it follows that $\cos \alpha(r) \in [1/2,1]$, therefore $\alpha(r) \in [0,\pi/3]$. Then a quick computation shows that
\[ |T_{ABC}| = \pi r^2+3r^2(\tan \alpha(r) -\alpha(r))=3r^2\left(\frac{\pi}{3}+\tan \alpha(r)-\alpha(r)\right)=:f(r).\]
Elementary computations show that
\[ \alpha(r)' = \frac{-1}{(1-r)^2 \sin \alpha(r)}, f'(r) = 6r \left(\frac{\pi}{3}-\alpha(r)\right)+3(2r\tan \alpha(r)-\sin \alpha(r)), \]
where $'$ denotes differentiation with respect to $r$. The terms appearing in the derivative $f'(r)$ are non-negative since $\alpha(r) \in [0,\pi/3]$ and
\[ 2r \tan \alpha(r)-\sin \alpha(r) = (1-2r)\sin \alpha(r)\geq 0.\]
Moreover, inequality is strict for $r \in (1/3,1/2)$ showing that $f$ is strictly increasing on $[1/3,1/2]$. 

We conclude that $r \mapsto h(T_{ABC})$ is strictly decreasing with respect to $r$, thus $h(T_{ABC})$ is maximal for $r = 1/3$, that is when $T_{ABC}$ coincides with the equilateral triangle of unit width. Thus Theorem \ref{thm:Cheeger-minw} is proved by direct comparison with \trc  \  sets. Indeed, if $K$ is a convex shape with minimal width equal to one and then $K$ contains a \trc \ set $T_{ABC}(r)$. Therefore
\[ h(K) \leq h(T_{ABC}(r))\leq h(T_{ABC}(1/3)),\]
where $T_{ABC}(1/3)$ is the equilateral triangle of unit width.

\section{Minimizing the torsional rigidity among shapes with minimal width}
\label{sec:lamk}

Proposition \ref{prop:reduction-caps} implies that knowing the behavior of the torsional rigidity on \trc{} sets may help decide whether the equilateral triangle maximizes this functional among shapes with minimal width equal to $1$. While the area and Cheeger constant of $T_{ABC}$ depends only on $r$, this is no longer the case for the torsional rigidity. In the following we divide the analysis into two steps: proving that $T(T_{ABC})$ is minimal when $ABC$ is equilateral for fixed inradius $r$ and secondly, proving that among equilateral \trc{} shapes the equilateral triangle is minimal. 

\subsection{Analysis at fixed inradius}

Fix $r \in (1/3,1/2)$. Let us first investigate the behavior of $T(T_{ABC})$ when the vertex $A$ moves radially on a circle of radius $1-r$. It is classical that if $\Omega$ is a general shape and $\theta$ is an infinitesimal perturbation of $\partial \Omega$ then the shape derivative of the Dirichlet energy is $T'(\Omega)(\theta) = \int_{\partial \Omega} (\partial_n u)^2 \theta\cdot n$, where $u$ solves \eqref{eq:torsion} and $\partial_n u = \nabla u \cdot n$ is the normal derivative of $u$. See \cite{henrot-pierre-english} or \cite{sokolowski-zolesio} for more technical details. 
 
Let $AS, AT$ be the tangents from $A$ to the incircle $D_r$. Consider the perturbation of $T_{ABC}$ where the vertex $A$ moves on the circle of radius $1-r$ along with the tangents from $A$ to $1-r$. The infinitesimal vector field $\theta_A$ which characterizes the movement of the boundary of $T_{ABC}$ is simply the rotation around the origin on segments $AS, AT$, $\theta_A$ being equal to zero everywhere else on $\partial T_{ABC}$. Assuming the rotation from $S$ towards $T$ it follows that for $X \in [AS]$ we have $\theta_A(X) \cdot n = -r \cdot XS$ and for $X \in [AT]$ we have $\theta_A(X)\cdot n = r\cdot XT$. Writing $X(t)$ as a uniform parametrization of the segments $[AT], [AS]$, we find that
\[ T'(\theta_A) = r\int_{[AT]} X(t)T(\partial_n u)^2 dt - r \int_{[AS]} X(t)S (\partial_n u)^2 dt.\]
Therefore, the following proposition holds:

\begin{prop}
	\label{prop:one-vertex-move}
	For every $t \in [0,AT]$ denote points $X(t)\in [AT]$, $Y(t) \in [AS]$ such that $TX(t)=SY(t)=t$. If for every $t \in [0,AT]$ we have
	\[ (\partial_n u)^2(X(t))<(\partial_n u)^2(Y(t))\]
	then assuming the cap from $A$ is not touching the caps from $B$ and $C$,  $A$ can be perturbed to strictly increase (or decrease) $T(T_{ABC})$, the torsional rigidity of $T_{ABC}$.
\end{prop} 

The proof follows at once, since under the hypotheses considered the derivative $T'(\theta_A)$ has a well definite sign, therefore there exists a perturbation strictly increasing or decreasing the torsional rigidity. The goal in the following will be to show that in non-equilateral configuration there exists a vertex for which inequality in Proposition \ref{prop:one-vertex-move} holds. Techniques used in the following are similar to those in \cite{FragalaVelichkov19} using symmetrizations and comparisons involving harmonic functions. 

\begin{figure}
	\centering
	\includegraphics[width=0.4\textwidth]{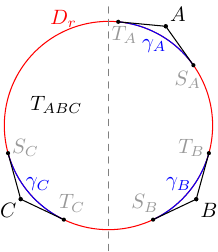}
	\caption{A three-$r$-cap shape $T_{ABC}$ together with the arcs $\gamma_A, \gamma_B,\gamma_C$ near points $A,B,C$. In the figure, when $AB<AC$, the symmetric of the left part with respect to the vertical axis is contained in the right part.}
	\label{fig:three-cap-arcs}
\end{figure}

Let us introduce the following notations:
\begin{itemize}[noitemsep]
	\item For vertex $A$ tangents to $D_r$ intersect the incircle $\partial D_r$ at points $T_A,S_A$ where $S_A$ comes first in the trigonometrical sense. If there is no confusion, the subscript is dropped. 
	\item The arc $(S_A,T_A)$ in $\partial D_r$ has angular length $2\arccos \frac{r}{1-r}$. Denote by $\gamma_A$ this arc, with uniform parametrization in the positive trigonometric sense $\gamma_A(t): [0,2\alpha] \to \partial D_r$. Arcs $\gamma_B$, $\gamma_C$ are defined in a similar manner. If needed, we use the notation $\gamma_A^{-}$ to refer to an arc parametrized in the opposite sense, still using uniform parametrization. See Figure \ref{fig:three-cap-arcs} for an illustration.
\end{itemize}

We are now ready to state the first comparison result in a non-equilateral configuration. 

\begin{prop}\label{prop:compBC}
	Assume $AB<AC$. Then for every $t \in [0,2\alpha]$ we have $u(\gamma_B(t))\geq u(\gamma_C(2\alpha-t))$. In other words, the torsion function $u$ takes larger values on $\gamma_B$ than on $\gamma_C^-$, point by point.
\end{prop}

\emph{Proof:} Assume, without loss of generality that $B,C$ are symmetric with respect to the $y$-axis, $x_B=-x_C>0$. Therefore $x_A>0$. Let $u>0$ be the solution of \eqref{eq:torsion} on $T_{ABC}$. Denote $K_- = T_{ABC} \cap \{x\leq 0\}, K_+ = T_{ABC} \cap \{x \geq 0\}$. We notice immediately that the symmetric of $K_-$ with respect to the $y$ axis is contained in $K_+$ since caps at $B$ and $C$ will overlap. Denote by $K$ the symmetric of $K_-$ with respect to the $y$ axis. 

For $(x,y) \in K$ define $\overline u(x,y) = u(-x,y)$ and let $w: K \to \Bbb{R}$ be given by $w = u-\overline u$. Then we have the following:
\begin{itemize}[noitemsep]
	\item $-\Delta w = 0$ since both $u$ and $\overline u$ verify $-\Delta u = -\Delta \overline u = 1$ in $K$. 
	\item $w = 0$ on $\partial K$ except, a non-trivial subset of the arc $\gamma_A$ where $w = u-\overline u$ will be strictly positive. 
\end{itemize}
Therefore, $w$ is a harmonic function with non-negative boundary conditions, strictly positive on a non-zero subset of $\partial K$. The strong maximum principle indicates that $w>0$ inside $K$. The result follows observing that $w>0$ on $\gamma_B$ and $w$ is the difference between values of $u$ on $\gamma_B$ and $\gamma_C^{-}$ on this arc. \hfill $\square$

The same type of comparison could indicate an inequality between gradients on segments arriving at $B$ and at $C$ in $T_{ABC}$. Notice that the result of the previous proposition does not allow the use of the shape derivative to conclude. Indeed, we need to compare the gradient of $u$ on adjacent segments. This is done in the following. The first step is a lemma regarding harmonic functions on the half disk with a sliding Dirichlet boundary condition. 

Denote by $D_+ = \{z \in \Bbb{C} : |z|\leq 1, \text{Im z} \geq 0\}$ the upper half disk in the complex plane. Consider $\delta>0$ and a function $\phi: [0,\delta] \to \Bbb{R}_+$ such that $\phi(0)=\phi(\delta) = 0$, $\phi(t)>0$ for $t \in (0,\delta)$. Consider $\theta \in (0,\pi-\delta)$. Denote by $\gamma_\theta$ the arc in the upper boundary $\Gamma= \{|z|=1, \text{Im z} \geq 0\}$ of $D_+$ corresponding to arguments in $[\theta, \theta+\delta]$. 

Consider the following harmonic function depending on $\theta$: 
\begin{equation}\label{eq:harmonic-slide} \begin{cases}
	\Delta u_\theta = 0 & \text{ in } D_+ \\
	u_\theta(e^{it}) = \phi(t-\theta)& \text{ on } \gamma_\theta \\
	u_\theta = 0 &\text{ on }\partial D_+\setminus \gamma_\theta.
\end{cases}
\end{equation}
Notice that the boundary condition on the arc $\gamma_\theta$ only depends on the fixed function $\phi$ defined on $[0,\delta]$. Also the arc $\gamma_\theta$ slides along the unit circle in the positive trigonometric sense as $\theta$ increases. We are interested in understanding the behavior of $u_\theta$ as $\theta$ increases, close to the upper boundary for points having small arguments.  

\begin{prop}
	\label{prop:harmonic-half-disk}
	Consider $\theta< \theta' \in (0,\pi-\delta)$ and let $u_\theta, u_{\theta'}$ be solutions of \eqref{eq:harmonic-slide} for corresponding parameters. For every $t \in (0,\theta)$ there exists $r_t \in (0,1)$ such that for $r \in (r_t,1)$ we have
	\[ u_{\theta'}(re^{it})<u_{\theta}(re^{it}).\]
\end{prop}

\begin{figure}
	\centering
	\includegraphics[width=0.7\textwidth]{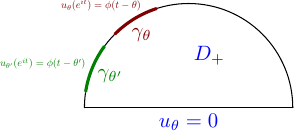}
	\caption{Illustration for boundary conditions for $u_\theta, u_{\theta'}$ given by \eqref{eq:harmonic-slide} for $\theta<\theta'$. Dirichlet boundary conditions $u_\theta=0$ and $u_{\theta'}=0$ hold everywhere on the boundary of the half-disk (including the diameter), except the arcs $\gamma_\gamma$ and $\gamma_{\theta'}$, respectively. On these arcs we have a boundary condition given by the same non-negative function $\phi$.}
	\label{fig:slide-disk}
\end{figure}

Intuitively, as the source term slides towards $\pi$, values of $u_\theta$ for points close to the unit circle having arguments smaller than $\theta$, decrease. See Figure \ref{fig:slide-disk} for an illustration. For the sake of clarity, the arcs $\gamma_\theta$ and $\gamma_{\theta'}$ do not intersect in the figure. However, the proof remains valid even when they intersect. 

\emph{Proof:} If $D$ is the unit disk then it is well known that for $f$ continuous on $\partial D$ (for simplicity) the harmonic function $u$ verifying $u=f $ on $\partial D$ is given by
\[ u(re^{it}) = \frac{1-r^2}{2\pi} \int_0^{2\pi} \frac{f(e^{is})}{|re^{it}-e^{is}|^2}ds = \frac{1-r^2}{2\pi} \int_{\partial D} \frac{f(z)}{|z-re^{it}|^2}ds.\]
See \cite[Chapter 6]{StraussPDE} for example. 
To find the explicit formula for the Dirichlet problem on $D_+$, consider an antisymmetric boundary data $f$ with respect to the imaginary part. It follows that for $r \in (0,1)$ and $r \in [0,\pi]$ we have
\begin{equation}
	\label{eq:kernel-half-disk}
	  u(re^{it}) = \frac{1-r^2}{\pi} \int_0^{2\pi}f(e^{is})\left(  \frac{1}{|re^{it}-e^{is}|^2}-\frac{1}{|re^{-it}-e^{is}|^2}\right)ds.
\end{equation}
Denote the integral kernel by 
\[ Q(r,t,s)=  \frac{1}{|re^{it}-e^{is}|^2}-\frac{1}{|re^{-it}-e^{is}|^2}.\]
Observing that the kernel corresponds to $\frac{1}{|z-x|^2}-\frac{1}{|z-\overline x|^2}$ where $z \in \partial D_+$ and $x \in D_+$ we find that $Q(r,t,s)\geq 0$. Moreover, since the distance between $x=re^{it}$ and $\partial D_+$ goes to $0$ as $r \to 1$ it follows that $\max_{s \in [0,\pi]} Q(r,t,s) \to +\infty$ when $t$ is fixed and $r \to 1$. 

A straightforward computation shows that
\[ Q(r,t,s) = \frac{1}{1+r^2-2r\cos(s-t)}-\frac{1}{1+r^2-2r\cos(s+t)}.\]
For fixed $r$ and $t$, denoting $g(s) = 2rQ(r,t,s)$, we have, denoting $C = \frac{1+r^2}{2r} \geq 1$, the following:
\[ g(s) = \frac{1}{C-\cos s \cos t-\sin s \sin t}-\frac{1}{C-\cos s \cos t+\sin s \sin t}.\]
The derivative is
\[ g'(s) = -\frac{\sin s\cos t-\cos s \sin t}{(C-\cos s \cos t-\sin s \sin t)^2}+\frac{\sin s \cos t+\cos s \sin t}{(C-\cos s \cos t+\sin s \sin t)^2}.\]
 The critical points of $g$ in $[t,\pi]$ verifying $g'(s)=0$ verify the equation
\[ \cos s(1+\cos^2 t+C^2)-\cos^3 s = 2C\cos t.\]

Replacing $\cos s$ with $X \in [-1,1]$, consider the function
\[ h : X\mapsto  X^3+2C\cos t- (1+\cos^2 t+C^2)X .\]
A classical one dimensional variational study will show that $h$ has a unique root $X_0 \in [-1,1]$. 
\begin{itemize}[noitemsep]
	\item  $h'(X) = 3X^2-1-\cos^2 t - C^2$ and $h''(X) = 6X$. Therefore $h'$ is strictly decreasing on $[-1,0]$ and strictly increasing on $[0,\cos t]$. It follows that, $h$ is strictly concave on $[-1,0]$ and strictly convex on $[0,1]$.
	\item Computations show that 
	\[ h(-1) = (C+\cos t)^2>0,h(1) = -(C-\cos t)^2<0,\]
	therefore, $h$ can only have one root in $[-1,1]$. 
\end{itemize}

Therefore $h$ vanishes exactly once in $[-1,1]$, let $X_0$ be the unique solution. Previous observations show that $h(X)>0$ on $[-1,X_0)$ and $h(X)<0$ on $(X_0,1]$. Moreover, 
\[ h(\cos t) = -\cos t(C-1)^2.\]
Therefore $X_0 \leq \cos t$ if $t\in (0,\pi/2]$ and $X_0 \geq \cos t$ if $t \in [\pi/2, \pi]$.

 Then $g'(s)=0$ is equivalent to $\cos s = X_0$ and since $s \in [0,\pi]$, the unique critical point $s_0 \in [0,\pi]$ is determined uniquely. We can see directly that $g'(t)>0$ and $g'(\pi)<0$. Therefore $g$ is increasing on $[0,s_0]$ and decreasing on $[s_0,\pi]$.
 
 Moreover, if $t \leq \pi/2$ then $s_0 \geq t$ and if $t \geq \pi/2$ we have $s_0 \leq t$. Nonetheless, since $C = (1+r^2)/(2r)$ we have that $r\to 1$ implies $C \to 1$. Since $h(\cos t) = -\cos t (C-1)^2$ we find that when $r \to 1$ the critical point $s_0$ of $g$ converges to $t$. 

Summarizing the previous observations we obtained that for fixed $r,t$ determining $x = re^{-it} \in D_+$ we have:
\begin{itemize}[noitemsep]
	\item $s \mapsto Q(r,t,s)$ is strictly decreasing on $[s_0,\pi]$, $s_0 > t$ if $t<\pi/2$ and $s_0 \to t$ when $r \to 1$.
	\item $s \mapsto Q(r,t,s)$ is strictly decreasing on $[s_0,\pi]$ and $s_0\leq t$ when $t \geq \pi/2$.
\end{itemize}

Consider now $t<\theta<\theta'$. We have, using the representation formula,
\[ u_\theta(re^{it})=\frac{1-r^2}{2\pi} \int_\theta^{\theta+\delta} \phi(s-\theta) Q(r,t,s)ds.\]
For fixed $t$, take $r'\in (0,1)$ large enough such that the critical point $s_0$ of $s \mapsto Q(r,t,s)$ verifies $s_0<\theta$. This is always possible, since $t<\theta$ and $s_0$ converges to $t$ as $r \to 1$. For every $r>r'$ we know that $s \mapsto Q(r,t,s)$ is strictly decreasing on $[\theta,\pi]$. Therefore
\[ u_\theta(re^{it})=\frac{1-r^2}{2\pi} \int_\theta^{\theta+\delta} \phi(s-\theta) Q(r,t,s)ds>\frac{1-r^2}{2\pi} \int_{\theta'}^{\theta'+\delta} \phi(s-\theta') Q(r,t,s)ds= u_{\theta'}(re^{it}).\]
\hfill $\square$ 

Building on the previous result, it is straightforward to arrive at the same conclusion, assuming $u_\theta$ and $u_{\theta'}$ verify, in addition another common Dirichlet boundary condition of the type
\[ u_\theta(e^{is}) = u_{\theta'}(e^{is}) = \psi(s), \text{ for every } s \in (0,t),\]
where $\psi: [0,t]\to \Bbb{R}$ is an arbitrary one dimensional function. 

\begin{prop}\label{prop:harmonic-half-disk-extra}
	Let $u_\theta$ be solution of 
	\begin{equation}\label{eq:harmonic-slide-modif} \begin{cases}
			\Delta u_\theta = 0 & \text{ in } D_+ \\
			u_\theta(e^{is}) = \phi(s-\theta)& s \in [\theta,\theta+\delta] \\
			u_\theta(e^{is}) = \psi(s) & s \in [0,t] \\
			u_\theta = 0 &\text{ otherwise on  }\partial D_+.
		\end{cases}
	\end{equation}
	where $\phi$ verifies conditions in Propostion \ref{prop:harmonic-half-disk} and $\psi$ is arbitrary.
	
	Then for every $t <\theta \in (0,\pi)$ there exists $r_t \in (0,1)$ such that for $r \in (r_t,1)$ and $\theta'>\theta$ we have
	\[ u_{\theta'}(re^{it})<u_{\theta}(re^{it}).\]
\end{prop}

\emph{Proof:} The proof follows the same lines as Proposition \ref{prop:harmonic-half-disk}, using the representation formula \eqref{eq:kernel-half-disk}, which contains a common part corresponding to $z$ on $\partial D_+$ having arguments in $[0,t]$.
\hfill $\square$

\begin{figure}
	\centering
	\includegraphics[width=0.7\textwidth]{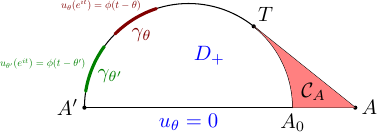}
	\caption{Geometric configuration for problem \eqref{eq:harmonic-K} on the domain $K$, the convex hull of the half disk and a point $A$ on the real axis. The boundary condition is $u_\theta=0$ except an arc $\gamma_\theta$ which slides towards $A'$. Values of $u_\theta$ decrease in $\mathcal C_A$ as the positive boundary condition slides.}
	\label{fig:slide-K}
\end{figure}

We are now ready to prove a similar result but on a modified domain which will help to answer the questions relevant to our main problem. Let $K$ be a modification of the upper half disk $D_+$ defined as follows:
\begin{itemize}[noitemsep]
	\item take $A$ on the real axis associated to the complex number $z_A>1$. Construct $AT$ the tangent from $A$ to $D_+$ with $T \in \partial D_+$. 
	\item Define $K$ as the convex hull of $A$ and $D_+$. 
	\item Denoting $A'$ the point associated to the complex number $z_{A'}=-1$ the boundary of $K$ is made of: the segment $AT$, the arc $\arc{TA'}$ on $\partial D_+$, the segment $A'A$ included in the real axis. 
	\item let $t_0 \in (0,\pi/2)$ such that $T = e^{it_0}$. Denote by $A_0$ (corresponding to $z=1$) the intersection of $AA'$ with the unit circle. 
	\item Denote by $\C_A$ the region of $K$ bounded by $AT$, the arc $\arc{A_0T}$ of $\partial D_+$ and $AA_0$. It is half of the cap determined by $A$ and $D_+$. See Figure \ref{fig:slide-K} for an illustration. 
	\item Like before, for $\delta>0$ and $\theta \in [t_0,\pi-\delta]$  consider the arc $\gamma_\theta \subset  \arc{TA'}$ associated to arguments in $[\theta,\theta+\delta]$ and an arbitrary positive function $\phi:[0,\delta] \to \Bbb{R}$, vanishing at $0$ and $\delta$.
\end{itemize}

Let $u_\theta$ be the solution to the problem
\begin{equation}\label{eq:harmonic-K}
	\begin{cases}
		-\Delta u_\theta = 0 & \text{ in }K\\
		u_\theta(e^{it}) = \phi(t-\theta) & \text{ on } \gamma_\theta \\
		u_\theta = 0 & \text{ otherwise on } \partial K.
	\end{cases}
\end{equation}
See Figure \ref{fig:slide-K} for an illustration. 

\begin{prop}\label{prop:harmonic-modif}
	Let $z$ be a point in the interior of $\mathcal C_A$. Then if $\theta<\theta'\leq \pi-\delta$ we have $u_\theta(z)>u_{\theta'}(z)$. As a consequence, the squared outward normal derivative of $u_\theta$ is larger than the squared outward normal derivative of $u_{\theta'}$ on the segment $AT$. 
\end{prop}

\emph{Proof:} Consider $\theta <\theta'\leq \pi-\delta$ and the notations above. Denote by $\psi$ the trace of $u_\theta$ on the arc $\arc{A_0T}$ of $\partial D_+$ corresponding to arguments in $[0,t_0]$. Denote by $w_\theta$ the solution of \eqref{eq:harmonic-slide-modif} with boundary data $\psi$ on $\arc{A_0T}$. Proposition \ref{prop:harmonic-half-disk-extra} shows that $w_{\theta'}<w_\theta$ close enough to the arc $\arc{A_0T}$ inside $D_+$. Hopf's boundary lemma shows that the normal derivatives across $\arc{A_0T}$ of $w_{\theta'}$ and $w_\theta$ verify the inequality $\partial_n w_{\theta'} >\partial_n w_\theta=\partial_n u_\theta$, assuming the normal pointing outwards on $\partial D_+$.

Let $v$ be the function defined on $K$ by 
\[ v = \begin{cases}
	w_{\theta'} & \text{ in } D_+ \\
	u_\theta & \text{ in } \C_A.
\end{cases}\]
where $u_\theta$ solves \eqref{eq:harmonic-K}. The function $v$ is piece-wise harmonic in $D_+$ and $\C_A$, continuous across the common boundary $\arc{A_0T} = D_+\cap \C_A$, with positive jump of the normal derivative across $\arc{A_0T}$.

Let $U \subset K$ be an arbitrary open set. Then for any non-negative $\phi \in C_c(U)$ we have $\int_U v \Delta \phi \leq 0$. Indeed, if $U\subset K$ or $U \subset \C_A$ then $v$ is harmonic and we have $\int_U v\Delta \phi = 0$. Otherwise, denote by $U_1 = U\cap D_+$, $U_2 = U\cap \C_A$ and $\Gamma = U \cap \arc{A_0T}$. Applying Green's formulas we find
\[ \int_U v\Delta \phi = \int_{U_1} v\Delta \phi +\int_{U_2} v\Delta \phi= \sum_{i=1,2}\int_{\partial U_i} v\frac{\partial \phi}{\partial n} - \phi \frac{\partial v}{\partial n}. \]
Since $\phi$ has compact support in $U$, only boundary integrals on the common boundary $\Gamma$ remain. Furthermore, $\phi$ is smooth across $\Gamma$ and $v$ is continuous. Therefore, since $\phi \geq 0$, 
\[ \int_U v\Delta \phi = -\int_{\Gamma} \phi (\partial_n w_{\theta'} - \partial_n u_\theta) \leq 0.\]
This implies that $v$ is superharmonic in $K$, having same Dirichlet boundary conditions as the harmonic function $u_{\theta'}$ on $\partial K$. Therefore $v \geq u_{\theta'}$ on $K$. 

In particular for $z \in \C_A$ we have $u_\theta(z) = v(z)\geq u_{\theta'}(z)$. The inequality is strict in a neighborhood of $\arc{A_0T}$ due to the inequality of normal derivatives. Since $u_\theta$ and $u_{\theta'}$ are harmonic in $\C_A$ the inequality is strict everywhere in the interior of $\C_A$. \hfill $\square$

\begin{prop}\label{prop:cap-comparison}
	Let $T_{ABC}$ be a \trc{} such that $AB<AC$. Then $T(T_{ABC})$ can be decreased by moving $A$ towards $C$ and $T_{ABC}$ does not minimize the torsional rigidity.
\end{prop}

\begin{figure}
	\centering
	\includegraphics[width=0.5\textwidth]{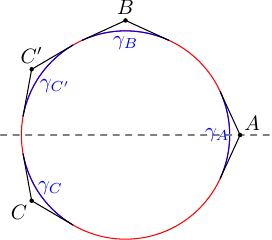}
	\caption{A three-$r$-cap shape with $AB<AC$ such that $A$ is on the positive real axis. The vertex $C$ is symmetrized with respect to the $x$-axis. The difference of the torsion function and its symmetrized will be compared in $\text{conv}(D_+,A)$.} 
	\label{fig:cap-comparison}
\end{figure}

\emph{Proof:} The proof is a simple consequence of Proposition \ref{prop:harmonic-modif}. The proposition was stated for the unit disk, but rescaling to disks of arbitrary radii preserves the conclusions. Assume that the incircle of $T_{ABC}$ is at the origin and $A$ lies on the positive real axis. An illustration of the geometric configuration is shown in Figure \ref{fig:cap-comparison}

\bo{Case 1.} Assume that the cap $\C_C$ is contained in the lower half plane $\{\text{Im}z \leq 0\}$. Since $AB<AC$, obviously the cap $\C_B$ is contained in the upper half plane. Recall the notations from Proposition \ref{prop:compBC}: $\gamma_B$ is the arc determined by tangents from $B$ to $D_+$. 

Reflect $C$ with respect to the real axis, obtaining a point $C'$. Denote by $\gamma_{C'}$ the arc on $\partial D_+$ determined by tangents from $C'$ to $D_+$. By construction, $\gamma_B$ and $\gamma_{C'}$ have same lengths. 

Let $u$ be the solution of the torsion equation \eqref{eq:torsion} on $T_{ABC}$. Let $K$ be the convex hull of the upper half disk of $D_r$ and $A$. For $x+iy \in K$ denote $\overline u(x,y) = u(x,-y)$. The difference $w = u-\overline u$ is harmonic in $K$. Consider the following harmonic functions:
\[ \begin{cases}
	-\Delta w_1 = 0 & \text{ on }K \\
	w_1 = u & \text{ on } \gamma_B \\
	w_1 = 0  & \text{ on } \partial K \setminus \gamma_B
\end{cases} 
\]

\[ \begin{cases}
	-\Delta w_2 = 0 & \text{ on }K \\
	w_2(re^{it}) = \overline u(re^{i(t+\sigma)}) & \text{ on } \gamma_B \\
	w_2 = 0  & \text{ on } \partial K \setminus \gamma_B
\end{cases} 
\]
where $\sigma$ is the angular difference between $\gamma_B$ and $\gamma_{C'}$,
\[ \begin{cases}
	-\Delta w_3 = 0 & \text{ on }K \\
	w_3(re^{it}) = \overline u & \text{ on } \gamma_{C'} \\
	w_3 = 0  & \text{ on } \partial K \setminus \gamma_{C'}
\end{cases}. 
\]
We have the following inequalities:
\begin{itemize}[noitemsep]
	\item $w_1 \geq w_2$ in $K$: Proposition \ref{prop:compBC} implies that values of on $\gamma_B$ are larger than values on $\gamma_C$ (with orientation reversed). Thus the boundary condition for $w_2$ is smaller than the boundary condition for $w_1$, implying the conclusion.
	\item $w_2 > w_3$ on $\C_A$: Consequence of Proposition \ref{prop:harmonic-modif}.
\end{itemize}

We conclude that $w_1 > w_3$ on $\C_A$ and therefore the Hopf boundary lemma implies that $\partial_n w_1>\partial_n w_3$ on the segment $TA$. Since $w = w_1-w_3 = u-\overline u$, the analogue inequality between normal derivatives of $u$ on tangent segments from $A$ to $D_r$ follows. Proposition \ref{prop:one-vertex-move} implies that $A$ is not in optimal position to minimize $T(T_{ABC})$.

\bo{Case 2.} $\C_C$ is not contained in the lower half plane $H_-=\{\text{Im} z \leq 0\}$, but $C \in H_-$. In this case symmetrizing across the real axis produces a less restricting situation. The arc $\gamma_{C'}$ is not completely contained in the upper half space $H_+$. The boundary condition for $w_3$ defined above will be replaced with a function smaller than $\overline u$ on $\gamma_{C'}\cap H_+$. The inequality between $w_2$ and $w_3$ becomes even stronger since $w_3$ has a pointwise smaller boundary condition than in the previous case. 

\bo{Case 3.} $C$ is contained in the upper halfspace $H_+ = \{\text{Im} z \geq 0\}$. In this case, the comparison with the symmetrized function of $u$ on $K$ is immediate since $u-\overline u$ is harmonic in $K$ with non-negative Dirichlet boundary conditions. \hfill $\square$

Previous results allow us to find the following.

\begin{thm}\label{thm:optimality-caps}
	Let $r \in (1/3,1/2)$ fixed. Among \trc{} sets $T_{ABC}$ with inradius $r$ the situation where $ABC$ is an equilateral triangle minimizes the torsional rigidity.
\end{thm}

\emph{Proof:} Note that $T_{ABC}$ is contained between two disks of fixed radii, therefore finite lower and upper bounds for the torsional rigidity $T(T_{ABC})$ can be given. The existence of a minimizer is assured by the Blaschke selection theorem for convex shapes and continuity properties of the torsion function for the Hausdorff convergence of convex shapes. See \cite[Chapters 2,4]{henrot-pierre-english}, \cite{henroteigs}, or \cite{AntunesBogosel22} for more details regarding existence theory.

Proposition \ref{prop:cap-comparison} shows that non-equilateral configurations are not maximal. Therefore the equilateral configuration is the maximal one. \hfill $\square$

\begin{rem}
	Incidentally, existence of a maximizing configuration and Proposition \ref{prop:cap-comparison} implies that the \trc{} $T_{ABC}$ which maximizes the torsional rigidity at fixed inradius $r \in (1/3,1/2)$ must consist of three adjacent caps $\C_A,\C_B,\C_C$.
\end{rem}

\subsection{Analysis of equilateral configurations} The previous sections imply that it is enough to investigate equilateral three cap configurations for inradius $r \in [1/3,1/2]$ to find the answer to Conjecture \ref{conj:torsion-minw}. Given $r \in [1/3,1/2]$ consider $ABC$ an equilateral triangle with circumradius $1-r$. The equilateral three-cap-set $T_{ABC}$ is the convex hull the triangle $ABC$ and of the disk of radius $r$ having the same center.  

\begin{figure}
	\centering
	\includegraphics[width=0.55\textwidth]{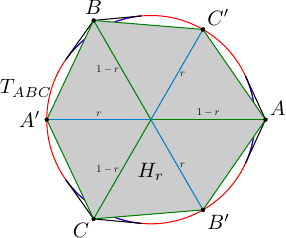}
	\caption{An equilateral \trc{} shape $T_{ABC}$ containing a hexagon $H_r$ with vertices $A,B,C,A',B',C'$ at alternating distances $r$ and $1-r$ from the center $O$.}
	\label{fig:equi-caps}
\end{figure}

Let $A',B',C'$ be the second intersections of $AO,BO,CO$ with the boundary of $T_{ABC}$. Denote by $H_r$ the convex hexagon obtained by taking the convex hull of $A,B,C,A',B',C'$. By construction $H_r$ is convex and contained in $T_{ABC}$. Moreover, the vertices of $H_r$ are at alternating distances $r,1-r$ from $O$. See Figure \ref{fig:equi-caps} for an illustration.

The monotonicity of the torsional rigidity shows that
\[ T(H_r) \leq T(T_{ABC}).\]
If in addition, the function $[1/3,1/2]\ni r \mapsto T(H_r)$ is strictly increasing, we will find that the torsional rigidity for the equilateral triangle is the smallest among all \trc{} sets, proving Theorem \ref{conj:torsion-minw}.

The uniqueness of the torsion function 
\[ -\Delta w = 1 , w \in H_0^1(H_r)\]
implies that $w$ is symmetric with respect to $AA', BB', CC'$. Thus, $w$ can be decomposed into $6$ components corresponding to the following Dirichlet-Neumann problems. 

For $r \in [1/3,1/2]$, let $XYZ$ be a triangle such that
\[ XY = r, XZ = 1-r.\]
and the angle at $X$ is $\angle YXZ = \pi/3$.
In the following $\Delta XYZ$ denotes the triangle $XYZ$. 
Consider $w$ the solution of
\begin{equation}\label{eq:torsion-tri}
	\begin{cases}
-\Delta w = 1   & \text{ in } \Delta XYZ \\
w = 0 & \text{ on }YZ \\
\partial_n w = 0 & \text{ on } XY \cup XZ.
	\end{cases}
\end{equation}
Symmetrizations of solutions of \eqref{eq:torsion-tri} recover the torsion function on $H_r$ which is a convex shape. Therefore $T(H_r) = 6\int_{\Delta XYZ} |\nabla w|^2$.

\begin{figure}
	\centering 
	\includegraphics[width=0.35\textwidth]{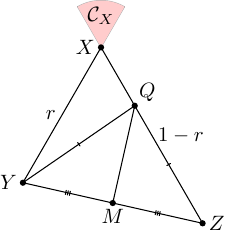}
	\caption{The triangle $XYZ$, one sixth of the hexagon $H_r$ illustrated in Figure \ref{fig:equi-caps}, having $\angle X = \pi/3$, $XY=r, XZ=1-r$, $r \in [1/3,1/2]$. The normal derivative of the torsion function is positive along $YQ$. Symmetrization arguments in Proposition \ref{prop:torsion-tri} show that the norm of the gradient of the torsion function $|\nabla w|$ in \eqref{eq:torsion-tri} is larger on $YM$ than on $MZ$.}
	\label{fig:triangleXYZ}
\end{figure}

We have, thus, the following properties:
\begin{itemize}[noitemsep,topsep=0pt]
	\item The maximum of $w$ is attained at the vertex $X$, for symmetry reasons and uniqueness of the maximizer of the torsion function on the convex domain $H_r$ \cite{Korevaar1987}. 
	\item The function $w^{1/2}$ is strictly concave on $\Delta XYZ$, with maximum at vertex $X$, therefore $w$ is strictly increasing on every segment of the form $X'X$ with $X' \in YZ$ \cite{Kennington1985}, \cite{Korevaar1987}. 
	\item The level lines of the function $w$ are convex and touch $XY, XZ$ normally. Therefore, for every $x \in \Delta XYZ$ we have $\nabla w(x) \in \C_X:= \{s(X-Y)+t(X-Z) : s,t \in \Bbb{R}_+\}$. 
	\item It follows that if $Q \in (XZ)$ then $\partial_n w>0$ where the normal is considered outward to $YQ$ in $\Delta YZQ$. This follows from the fact that $\C_X$ is contained in the halfplane determined by $YQ$ for which the scalar product with the normal is positive. See Figure \ref{fig:triangleXYZ}.
\end{itemize}

We are now ready to show how the torsion of $XYZ$ varies with $r$. 
\begin{prop}\label{prop:torsion-tri}
	Let $XYZ$ be a triangle with the properties described above, assuming $XY<XZ$, and let $Q \in (XZ)$ such that $YQ=ZQ$. Consider $M$ the midpoint of $YZ$.
	
	(a) If $U \in (MY), V \in (MZ)$ are such that $UM = VM$ ($U,V$ are symmetric with respect to the midpoint $M$) then 
	\begin{equation}\label{eq:ineq-sym}
		 (\partial_n w)^2 (U)> (\partial_n w)^2 (V).
	\end{equation}
	
	(b) The function $r \mapsto T(H_r)$ is strictly increasing on $[1/3,1/2]$.
\end{prop}

The geometric configuration is illustrated in Figure \ref{fig:triangleXYZ}.

\emph{Proof:} (a) Take a coordinate system such that $Y,Z$ are on the $x$-axis and $M$ is at the origin. Assume $Y$ is on the positive axis. Consider $\overline w$ defined on $\Delta MQY$ by symmetry, $\overline w(x,y) = w(-x,y)$. 

Recall that the properties of the torsion function on convex domains recalled above imply that $\partial_n w>0$ on $YQ$. Then $u = w-\overline w$ verifies the following conditions:
\begin{equation}
    \label{eq:tri-sym}
    \begin{cases}
    -\Delta u = 0 & \text{ on }\Delta MQY \\
     u   = 0 & \text{ on } MY \cup MQ \\
     \partial_n u >0 & \text{ on } YQ.
\end{cases}.
\end{equation}
It follows that $u>0$ in the interior of $\Delta MQY$. Hopf's boundary lemma shows that $\partial_n u<0$ on $MY$, which implies $\partial_n w<\partial_n \overline w$ on $MY$. Since the function $w$ is positive, the outward normal derivatives on $MY$ are negative. Therefore \eqref{eq:ineq-sym} follows. 

(b) Increasing $r$ only produces normal boundary perturbations on the segment $YZ$. More precisely: if $r \gets r+\delta r$ then
\begin{itemize}[noitemsep]
	\item the normal infintesimal displacement at $Y$ is given by $\cos \angle(\overrightarrow{ XY},n_{YZ})=\sin \angle XYZ$.
	\item the normal infinitesimal displacement at $Z$ is given by $-\cos \angle(\overrightarrow{ XZ},n_{YZ})=-\sin \angle XZY$
\end{itemize}

Therefore, since $XY<YZ$ the normal displacement at $Y$ is larger in absolute value than the one at $Z$, hence there exists a point $M' \in (MZ)$ such that the induced infinitesimal normal boundary perturbation $\theta\cdot n$ is strictly positive on $YM'$ and strictly negative on $M'Z$.

The shape derivative of the torsion with respect to movement on the Dirichlet parts of the boundary is still given by $T'(\Delta XYZ)(\theta) = \int_{YZ} (\partial_n w)^2 \theta\cdot n$. Inequality \eqref{eq:ineq-sym} and the discussion above shows that $|(\partial_n w)^2 \theta \cdot n |$ is larger on $YM$ than on $ZM$ for points symmetric with respect to $M$. This implies that $T'(\Delta XYZ)(\theta)>0$ for perturbations associated to increasing the parameter $r$. \hfill $\square$

The proof of the main result follows combining the arguments above.

\noindent\bo{Proof of Theorem \ref{conj:torsion-minw}}:

 Monotonicity of the torsional rigidity allows to reduce the problem to \trc{} shapes, according to Proposition \ref{prop:reduction-caps}. Theorem \ref{thm:optimality-caps} shows that among \trc{} shapes with fixed inradius, the equilateral one minimizes the torsion. 

Each equilateral \trc{} shape contains a hexagon $H_r$ defined above (see Figure \ref{fig:equi-caps}). Moreover, the torsion of these hexagons is minimal for $r = 1/3$ as shown in Proposition \ref{prop:torsion-tri}. Finally, by direct comparison, the equilateral triangle minimizes the torsion among shapes with given minimal width. \hfill $\square$

\begin{rem}
	The optimality of the equilateral triangle for the maximization of the first Dirichlet-Laplace eigenvalue remains open, although numerical simulations in \cite{Bogosel-Convex} and the results of this paper, regarding the minimization of the torsional rigidity are strong indicators regarding its validity. 
	
	Numerical simulations show that Proposition \ref{prop:cap-comparison} and Proposition \ref{prop:torsion-tri} remain valid for Dirichlet-Laplace eigenvalues (with the reversed monotonicity considerations). Therefore, the same proof strategy could apply. Results of \cite{FragalaVelichkov19} also motivate this assumption.
	
	From a theoretical point of view, Proposition \ref{prop:compBC} can be proved for the Dirichlet-Laplace eigenvalue using similar comparison arguments using partial symmetrization. Similarly, the proof of Proposition \ref{prop:torsion-tri} could be adapted for the eigenvalues. The analysis involving harmonic functions like in Theorem \ref{eq:harmonic-slide-modif} does not extend directly to the case of eigenvalues.
\end{rem}

\section{Conclusions}

In this paper, we show how two shape optimization problems under minimal width constraint can be solved by direct comparison with particular shapes. Each minimal width shape contains a \trc{} shape, having a simple structure, depending on three parameters. Restricting the study to these shapes allows to further reduce the problem to equilateral configurations, yielding a one dimensional problem which is then solved.

We show that the Cheeger constant is maximized by the equilateral triangle, recovering a result from \cite{Cheeger-minw}. We also prove that the torsional rigidity is minimized by equilateral triangle under minimal width constraint. This is a first case where a shape optimization problem involving PDEs under a minimal width constraint is solved analytically. 

The techniques involved in the proof are mainly comparison results for harmonic functions, in the same spirit as \cite{FragalaVelichkov19}. It is likely that similar considerations might lead to a positive result to Conjecture \ref{conj:lamk-minw}, regarding the maximization of the first Dirichlet-Laplace eigenvalue under width constraint.

\bigskip

\noindent \bo{Data availability statement:} No data was generated or used in this article. 

\noindent \bo{Conflicts of interest:} The authors have no relevant financial or non-financial interests to disclose.

\noindent \bo{Acknowledgement:} The author was partially supported by the ANR Project:
STOIQUES. 

\bibliographystyle{abbrv}
\bibliography{../DiamConstr/biblio}

\bigskip
\small\noindent
Beniamin \textsc{Bogosel}

\noindent Centre de Math\'ematiques Appliqu\'ees, CNRS,\\
\'Ecole Polytechnique, Institut Polytechnique de Paris,\\
91120 Palaiseau, France \\
{\tt beniamin.bogosel@polytechnique.edu}\\
{\tt \nolinkurl{https://beniamin-bogosel.github.io/}}

\bigskip

\noindent Faculty of Exact Sciences,\\ 
Aurel Vlaicu University of Arad,\\
2 Elena Dr\u agoi Street,\\
 Arad, Romania\\
 {\tt beniamin.bogosel@uav.ro}\\

\end{document}